\documentclass[11pt]{amsart}
\usepackage{amssymb}
\usepackage{graphicx}
\renewcommand{\a}{\alpha}
\renewcommand{\b}{\beta}
\renewcommand{\d}{\delta}

\newcommand{\e}{\varepsilon}

\renewcommand{\l}{\lambda}
\newcommand{\s}{\sigma}

\renewcommand{\o}{\omega}
\newcommand{\p}{\phi}

\newcommand{\G}{\Gamma}

\renewcommand{\P}{\Phi}

\newcommand{\PP}{\mathbb{P}}
\newcommand{\Z}{\mathbb{Z}}
\newcommand{\C}{\mathbb{C}}

\newcommand{\CP}{\mathbb{CP}}

\newcommand{\iso}{\cong}

\newcommand{\im}{\operatorname{im}}

\newcommand{\too}{\longrightarrow}

\newcommand{\Diff}{\operatorname{Diff}}
\newcommand{\vol}{\operatorname{vol}}
\newcommand{\Map}{\operatorname{Map}}

\newcommand{\sref}{s^{\text{\scriptsize ref}}}
\newcommand{\bd}{\partial}
\newcommand{\bbd}{\bar{\partial}}
\newcommand{\x}{\times}
\newcommand{\ox}{\otimes}

\newtheorem{proposition}{Proposition}[section]
\newtheorem{theorem}[proposition]{Theorem}
\newtheorem{definition}[proposition]{Definition}
\newtheorem{lemma}[proposition]{Lemma}

\newtheorem{corollary}[proposition]{Corollary}
\newtheorem{remark}[proposition]{Remark}
\newtheorem{question}[proposition]{Question}

\title{Generic behavior of asymptotically holomorphic Lefschetz pencils.}

\thanks{\noindent First author supported by CICYT BFM2003--06001. \\
Second and third authors supported by CICYT BFM2000--0024.
\\ Third author supported by Fundaci\'on Pedro Barri\'e de la
Maza.
\\ Partially supported by The European Contract Human Potential
Programme, Research Training Network HPRN-CT-2000-00101.}

\date{April, 2004. Revised June 2004}

\author{Jaume Amor\'os}
\address{Departament de Matem\`{a}tica Aplicada I \\
Universitat Polit\`{e}cnica de Catalunya \\ 08028 Bar\-ce\-lo\-na \\
Spain} \email{jaume.amoros@upc.es}

\author{Vicente Mu\~noz}
\address{Departamento de Matem\'aticas \\
Universidad Aut\'onoma de Madrid
\\ 28049 Madrid \\ Spain}
\email{vicente.munoz@uam.es}

\author{Francisco Presas}
\address{Departamento de Matem\'aticas \\
Universidad Aut\'onoma de Madrid
\\ 28049 Madrid \\ Spain}
\email{francisco.presas@uam.es}

\begin{document}

\renewcommand{\theenumi}{\roman{enumi}}

\begin{abstract}
We prove that the vanishing spheres of the Lefschetz pencils
constructed by Donaldson on symplectic manifolds of any dimension
are conjugated under the action of the symplectomorphism group of
the fiber. More precisely, a number of generalized Dehn twists may
be used to conjugate those spheres.
 This implies the non-existence of homologically trivial
 vanishing spheres in these pencils. To develop the proof,
 we discuss some basic
 topological properties of the space of asymptotically holomorphic
 transverse sections.
\end{abstract}

\maketitle

\section{Introduction} \label{sec:introduction}
In this article  we analyze the generic behavior of vanishing
spheres in the symplectic pencils introduced by Donaldson in
\cite{Do99}, henceforth referred to as {\em Donaldson's
$\e$-transverse Lefschetz pencils\/} (see Section \ref{asym_tools}
for precise definitions), and show it to be similar to the case of
Lefschetz pencils for complex projective varieties. Using the
pencils as a tool, we start the study of the symplectic analogue
of the dual variety in algebraic geometry, which we believe will
be of interest in symplectic topology.

\medskip

The property of Lefschetz pencils on projective
algebraic varieties that we seek to extend is classical:

\begin{theorem} [cf. SGA 7 XVIII, 6.6.2]  \label{main:alg}
 Let $M$ be a complex projective manifold, and $L \to M$ an
 ample line bundle.
 For $k$ large enough the pencils associated to quotients of
 holomorphic sections of the line bundle $L^{\ox k}$ have vanishing
 spheres which are conjugated under the action of the group
 $\Diff^+(F)$ of orientation-preserving diffeomorphisms of the
 generic fiber $F$.
\end{theorem}

The main result of this article is a generalization of
the previous theorem to the symplectic case. By establishing the
irreducibility of symplectic analogues to the dual variety, defined
in \eqref{def:sect}, we will prove a more precise version, which replaces the
diffeomorphism group of $F$ with its smallest possible subgroup.
For that we need the
following

\begin{definition} Two Lagrangian submanifolds $L_0$ and $L_1$ are Lagrangian
isotopic if there is a continuous family of Lagrangian
submanifolds $\{ L_t\}_{t\in[0,1]}$ connecting them.
\end{definition}

Then, we show

\begin{theorem} \label{moser}
 For $k$ large enough Donaldson's $\e$-transverse Lefschetz pencils
 have vanishing spheres that, up to Lagrangian isotopy, are
 conjugated under the action of the group of symplectomorphisms
 generated by the Dehn twists of the fiber along the vanishing
spheres of the pencil itself.
\end{theorem}

We point out that the conjugating group of symplectomorphisms of
Theorem \ref{moser} is not elementary, as Donaldson's pencils
exist for large $k$ and have $O(k^{\dim M/2})$ singular fibers.
Nevertheless, this group may be substantially smaller than the
full symplectomorphism group of the fiber. Note that its action on
the homology of the fiber is trivial except in middle degree,
while in many instances, such as products and bundles of
symplectic manifolds, there exist symplectomorphisms acting
nontrivially in other cohomology groups.

Finally we will give the following application of theorem \ref{moser}:

\begin{theorem} \label{non-triv}
 For $k$ large enough Donaldson's $\e$-transverse
 Lefschetz pencils satisfy that all the vanishing spheres are
 homologically non-trivial.
\end{theorem}

In the case of dimension $4$, I. Smith \cite{Ivan} has already
proved that the vanishing spheres appearing in the Donaldson's
$\e$-transverse Lefschetz pencils are homologically non-trivial
for $k$ even. On the other hand, for general symplectic Lefschetz
pencils in $4$-manifolds \cite{BK} gives a lower bound for the
number of homologically non-trivial vanishing spheres.

In the cases of dimension greater than $4$, Theorem \ref{non-triv}
is still a meaningful result, as homologically trivial Lagrangian
spheres are known to exist in higher dimensional symplectic
manifolds at least in dimension $4k+2$ (see the examples of
\cite[3.d]{STY} and \cite{ALP}, for
instance), and there is great flexibility in constructing
symplectic Lefschetz pencils adapted to given Lagrangian
submanifolds (\cite{AMP}).

Theorem \ref{main:alg} is a consequence of the irreducibility of
the {\em dual
variety}, which is the subset of the dual projective space defined
by hyperplanes tangent to an embedded projective variety. It is
not too hard to give a definition of this variety in the
symplectic case, by combining ideas of \cite{MPS00} and
\cite{Do99}. This dual ``symplectic'' variety
appears as an ``asymptotically holomorphic'' divisor. The divisor
is not smooth and, in particular, has self-intersections. It is
difficult to control the behavior of
self-intersections for this kind of divisor (for an example
without known solution see \cite{AK00}). Therefore we have chosen
an alternative route to prove Theorem \ref{moser}.
However it makes sense to speak in a ``rough sense'' of the
complement of the dual variety, i.e., the space of transverse
hyperplanes. We will give some topological properties of this
space.

The organization of this paper is as follows. In Section
\ref{asym_tools} we give the definitions and results introduced in
\cite{Do96,Do99} needed to carry out the proofs.
In Section \ref{case:symp} we establish our symplectic analogue to the
irreducibility of the dual variety, on which the whole paper rests.
A discussion of the topology of the set of ``transverse'' sections
is carried out
in Section \ref{comp_dual}. Using this, Theorems \ref{moser} and
\ref{non-triv} are proved in Section \ref{case:moser}. Finally, in
Section \ref{invariant} we briefly discuss the possibility of obtaining a
symplectic invariant from the space of transverse sections.

\vskip5pt

{\bf Acknowledments:} Some of the results contained in Section
\ref{comp_dual} are the result of conversations between D. Auroux
and the third author. In particular, the statement of Theorem
\ref{thm:topo} was suggested by him to us. We thank Universidad
Aut\'onoma de Madrid for its hospitality to the first author
during the elaboration of this work. We are also grateful to L.
Katzarkov and I. Smith for useful comments.

\section{Donaldson's asymptotically holomorphic theory} \label{asym_tools}

Let $(M,\o)$ be a symplectic manifold of integer class, i.e.,\
satisfying that $[\o]/2\pi$ admits an integer lift to $H^2(M;\Z)$.
We define $L$ to be a hermitian bundle with connection whose
curvature is $R_L=-i\o$. Moreover fixing a compatible almost
complex structure $J$, this defines a Riemannian metric
$g(u,v)=\o(u,Jv)$ on $M$. We also consider the sequences of
metrics $g_k=kg$, $k\geq 1$. The following definitions are needed.

\begin{definition}
A sequence of sections $s_k$ of the bundles $L^{\ox k}$ has
$C^3$-bounds $c$ if it satisfies the following bounds \textup{(}in
$g_k$-metric\textup{)}
 $$
 |s_k|\leq c, \qquad |\nabla^r s_k|\leq c \quad
 (r=1,2,3), \qquad |\nabla^r \bbd s_k|\leq ck^{-1/2} \quad (r=1,2).
 $$
\end{definition}

\begin{definition}
A sequence of sections $s_k$ of the hermitian bundles $E_k$ is
$\e$-transverse to zero at the point $x\in M$ if one the
two following conditions is fulfilled
\begin{enumerate}
\item $|s_k(x)|\geq \e$.
\item $\nabla s_k:T_xM \to (E_k)_x$ is surjective and it admits
a right inverse of norm bounded by $\e^{-1}$, using the norm $g_k$ in $M$.
\end{enumerate}
\end{definition}

The sequence is said to be transverse to zero in a set $U$ if it
is transverse to zero at every point of the set.

We need to introduce notations to control $1$-parametric families
of sections. In particular $(J_t)_{t\in[0,1]}$ will denote a
$1$-parametric family of almost complex structures compatible with
$\o$. We can associate to this family a sequence of families of
metrics $g_{k,t}$. A family of sequences of sections $s_{k,t}$ has
$C^3$-bounds $c$ if $s_{k,t}$ has $C^3$-bounds $c$ with respect to
the almost complex structure $J_t$ and the metric $g_{k,t}$ for
each $t\in [0,1]$. (Note that such a family is a continuous path
in the space of sections of $L^{\ox k}$ with the $C^3$-topology.)
The same remark must be applied for families of $\e$-transverse
sequences.

There exist local objects that are close to being holomorphic, as
given by the following two lemmas

\begin{lemma} [Lemma 3 in \cite {Au00}] \label{lem:asymp_Darboux}
Let $x_t\in M$,  $t\in [0,1]$, be a path in $M$. Then there exist
complex Darboux coordinates depending continuously on $t$,
$\Phi_t: B_{g_k}(x_t, k^{1/2}) \to B_{\C^n}(0, k^{1/2})$,
$\Phi_t(x_t)=0$ such that the inverse $\Psi_t= \Phi_t^{-1}$ of the
coordinate map is nearly pseudo-holomorphic with respect to the
almost complex structure $J_t$ on $M$ and the canonical complex
structure on $\C^n$. Namely, the map $\Psi_t$ satisfies $|\nabla^r
\Psi_t|=O(1)$ for $r=1,2,3$ and $|\bbd \Psi_t(z)|=
O(k^{-1/2}|z|)$, and $|\nabla^r \bbd \Psi_t|=O(k^{-1/2})$ for
$r=1,2$.
\end{lemma}

We also have

\begin{lemma} [Lemma 3 in \cite{Au97}] \label{lem:exis_sect}
There exist constants $\l>0$ and $c_s>0$ such that, given any
continuous path $x_t : [0,1] \to M$, and large $k$, there exist
sections $\sref_{k,x_t}$ of $L^{\ox k}$ over $M$ with the following
properties:
\begin{enumerate}
\item The sections $s_{k,x_t}^{\text{\scriptsize ref}}$ have $C^3$-bounds
with respect to $J_t$, independent of $t$.
\item They depend continuously on $t$.
\item The bound $|\sref_{k,x_t}|>c_s$ holds over the ball
of $g_k$-radius $1$ around $x_t$.
\item $|s_{k,x_t}^{\text{\scriptsize ref}}(q)| \leq
\exp(-\l d_k(x_t,q)^2)$, where $d_k$ is the distance associated to
$g_k$.
\item $s_{k,x_t}^{\text{\scriptsize ref}}$ is supported in a ball
around $x_t$ of $g_k$-radius $c k^{1/6}$ for some constant $c$.
\end{enumerate}
\end{lemma}

\medskip

Recall from \cite{Do99}:
\begin{definition} A symplectic Lefschetz pencil
on $(M,\o)$ is a surjective map $\p:M-N \to \CP^1$, with $N$
a codimension 4 symplectic submanifold, such that every
point $p \in M$ has a complex--valued coordinate neighbourhood
$\psi: U\subset M \to \C^n$ sending $p$ to $(0,\dots ,0)$,
the standard symplectic structure $\omega_0$ of $\C^n$ to $\o$,
and such that
\begin{enumerate}
\item For $p \in N$, $N$ has local equation
 $\{z_1=z_2=0 \}$, and $\p(z_1, \ldots, z_n)=z_2 / z_1$.
\item For finitely many critical points $p_1, \dots, p_{\Lambda} \in M-N$,
$\p(z_1, \ldots, z_n)=z_1^2+
 \cdots + z_n^2$ (the {\em ordinary quadratic singularity} in algebraic
geometry).
\item For all other points $p \in M-N$, $p \neq p_1, \dots, p_{\Lambda}$,
 $\p(z_1, \ldots, z_n)= z_1$.
\end{enumerate}
\end{definition}

The main result of \cite{Do99} is

\begin{theorem} [Theorem 2 in \cite{Do99}] \label{thm:exist}
 Given a symplectic manifold $(M,\o)$ such that the cohomology
 class $[\o]/2\pi$ has an integer lift to $H^2(M,\Z)$,
 there exists a symplectic Lefschetz pencil whose fibers
 are homologous to the Poincar\'{e} dual of $k[\o]/2\pi$, for
 $k$ large enough.
\end{theorem}

Donaldson's proof of Theorem \ref{thm:exist} goes through the
following two steps:

\begin{proposition} \label{prop:get_tr}
Given a sequence of sections $s_{k}^1\oplus s_{k}^2$ of $L^{\ox
k}\bigoplus L^{\ox k}$ which has $C^3$-bounds $c$, then for any
$\d>0$, there exists $\e>0$ and a sequence of sections
$\s_{k}^1\oplus \s_{k}^2$ such that the following conditions are
satisfied for $k$ large
\begin{enumerate}
\item $|s_{k}^j-\s_{k}^j|_{C^1}\leq \d$, for $j=1,2$.
\item $\s_{k}^1$ is $\e$-transverse to zero over $M$.
\item $\s_{k}^1\oplus\s_{k}^2$ is $\e$-transverse to zero over $M$.
\item Denoting by $Z_{k,\e}=\{ p\in M: |\s_{k}^1|\leq \e \}$, the
map $\bd \left( \s_{k}^2 / \s_{k}^1\right)$ is $\e$-transverse to
zero in $M-Z_{k,\e}$.
\end{enumerate}
Moreover, given a $1$-parametric family of sequences of sections
$s_{k,t}^1\oplus s_{k,t}^2$, $t\in [0,1]$, of $L^{\ox k}\bigoplus
L^{\ox k}$ which have $C^3$-bounds $c$, then there exists a family
of sequences of sections $\s_{k,t}^1\oplus \s_{k,t}^2$ satisfying
the properties above for each $t\in [0,1]$. The result also holds
for continuous families of sequences of sections parametrized by
$t\in S^1$.
\end{proposition}

The second step in the proof is

\begin{proposition} \label{prop:get_model}
Given $\rho >0$ and a sequence of sections $s_{k}^1\oplus s_{k}^2$ satisfying
the last three properties of Proposition \ref{prop:get_tr}, then
the maps
 $$
 \p_{k}:M-Z(s_{k}^1\oplus s_{k}^2) \to \CP^1
 $$
are well defined and for $k$ large enough there exists a perturbation
producing a map
$\hat{\p}_{k}$ which defines a symplectic Lefschetz pencil
and verifies $\|\hat{\p}_{k}-\p_{k} \|_{C^1,g_k} <\rho$.

The same result holds for continuous families of sequences of
sections $s_{k,t}^1\oplus s_{k,t}^2$, $t\in [0,1]$ or $t\in S^1$,
satisfying the
conditions of Proposition \ref{prop:get_tr}.
\end{proposition}

For proofs of these results we refer the reader to \cite{Do99}.
The $1$-parametric version with parameter space $S^1$ follows
easily from the version with parameter space $[0,1]$.

We will call {\em Donaldson's
$\e$-transverse Lefschetz pencils} the Lefschetz pencils yielded
by Proposition \ref{prop:get_model}. The
$1$-parametric part of Proposition~\ref{prop:get_model} shows that
the symplectic Lefschetz pencils obtained by this procedure are
all isotopic for any given $k$ large enough.

\section{Vanishing spheres for symplectic Lefschetz pencils} \label{case:symp}

After blowing-up the symplectic submanifold $N$, any pencil
becomes a fibration $\tilde M \to \CP^1$. We can define a
canonical symplectic connection outside the critical points in
this fibration using the symplectic orthogonal of the tangent
space to the fiber as the horizontal subspace. The parallel
transport and geometric monodromy so defined act by
symplectomorphisms in the fibers (cf.\ \cite[Lemma 6.18]{MS94}).

It is possible to define a parallel transport ending in a singular
fiber, and the local computation of the Picard--Lefschetz formula
in algebraic manifolds (see \cite[vol.\ 2, ch.\ 1]{AVGZ}) shows
that in this case Lagrangian spheres in the regular fibers
contract to the singular point, forming a Lagrangian disk in the
manifold $M$. These spheres are called the {\em Lagrangian
vanishing spheres} in the regular fiber.

For a basepoint regular value $w\in \CP^1$, a Lagrangian vanishing
sphere in the fiber $F=p^{-1}(w)$ is determined by every path from
$w$ to a critical value $w_i$. The choice of a set of
non--intersecting paths $\gamma_1, \dots , \gamma_s$ from $w$ to
the singular values $w_1, \dots, w_s$ of the pencil defines a set
of vanishing spheres in $F$ that, together with a choice of isotopy
from the product of monodromies around all critical values to identity,
 determine the diffeomorphism type
of the pencil, in particular of the total space $\tilde M$
(\cite{Kas80}).

\medskip

Lefschetz pencils on a projective manifold $M \hookrightarrow \CP^N$
are defined by lines in the dual projective space $\CP^1 \cong l
\subset (\CP^N)^*$ intersecting the dual variety $M^*$ transversely
at smooth points. The hyperplane sections defined by the
line $l$ cover $M$, are smooth outside the finite intersection
$l \cap M^*$, where they consist of a hyperplane section with an
ordinary quadratic singularity (see SGA 7 XVII, XVIII for a complete
discussion).

The conjugacy modulo diffeomorphism of Theorem \ref{main:alg}
follows from the irreducibility of the dual variety $M^*$ for high
$k$: its smooth points form a connected open set $U \subset M^*$,
and given two critical values $w_i,w_j$ of the pencil, a path
connecting them in $U$ defines a smooth family of hyperplane
sections of $M$ whose parallel transport induces the conjugacy of
vanishing spheres.

In the symplectic case, instead of considering families of pencils
interpreted as lines in the dual projective space, we consider the
families in abstract. The result we want to prove is

\begin{proposition} \label{prop:hard}
The sequences of pencils provided by Theorem \ref{thm:exist}
satisfy, for $k$ large enough, that all Lagrangian vanishing
spheres of a generic fiber are equivalent under the action of an
orientation-preserving diffeomorphism of the fiber.
\end{proposition}

\begin{proof}
Choose a sequence of sections $s_k=s_k^1\oplus s_k^2$ satisfying
the conditions of Proposition \ref{prop:get_tr}, and let
$\p_k=s_k^2/s_k^1$ be the associated maps to $\CP^1$. For each $k$
select critical points $p_0$ and $p_1$ of the map $\bd
\p_k$. {}From \cite{Do99} we know that these points will become
the critical points of the pencil after the perturbation performed
in Proposition \ref{prop:get_model}. Now choose a path $p_t$ in
$M$ joining $p_0$ and $p_1$. We can choose this path to avoid a
$c'$-neighborhood of $Z(s_k^1)$ for some uniform $c'$, since by
\cite{Do99} the $\e$-transversality implies that $p_0$ and $p_1$
are away from a $c'$-neighborhood of $Z(s_k^1)$, for large $k$.
Also we may suppose that $p_t$ is stationary in $[0,\b]$ and
$[1-\b,1]$, for some small $\beta>0$. Also let $\l_t=\p_k(p_t)$,
which is a uniformly bounded path in $\C$.

Choose a family of sections $s_{k,p_t}^{\text{\scriptsize ref}}$
of $L^{\ox k}$ satisfying the properties of Lemma
\ref{lem:exis_sect} for the path $p_t$. We use these sections to
trivialize $L^{\ox k}$ in a neighborhood $U_t$ of fixed
$g_k$-radius $O(1)$ of $p_t$. The map $\p_k$ is given by
 \begin{eqnarray*}
  f_{k}: U_t & \to & \C \\
  q & \mapsto & \frac{s_{k}^2}{s_{k}^1}(q),
 \end{eqnarray*}
which has $C^3$-bounds $c_1 c$ (with the obvious adapted
definition) where $s_k$ has $C^3$-bounds $c$ and $c_1>0$ is a
constant depending only on the geometry of the manifold (not
depending on $k$, the initial sections, etc).

By Lemma \ref{lem:asymp_Darboux}, we may trivialize the
neighborhood $U_t$ by a chart $\P_{t} : B_{\C^n}(0,1) \to U_t$,
$\P_t(0)=p_t$, and
 $$
 \hat{f}_{k,t}=f_{k}\circ \Phi_{t}:B_{\C^n}(0,1) \to \C
 $$
has $C^3$-bounds $c_2 c$ in the unit ball of $\C^n$. Moreover
$\partial \hat{f}_{k,j}$ is $\e/c_2$-transverse to zero for
$j=0,1$. Again the constant $c_2>0$ depends only on the geometry
of the manifold $M$. Identifying the tangent space at the origin
to $\C^n$ we define $\hat{h}_{k,j}$, $j=0,1$,
as the quadratic form on $\C^n$ associated to $\partial \partial
\hat{f}_{k,j}(0)$. Moreover, it is possible to
construct a path of non-degenerate quadratic forms $\hat{h}_{k,t}$
starting and ending at the two previous quadratic forms. Also the
eigenvalues of the quadratic forms of the path can be bounded
below and above by the eigenvalues of the two initial quadratic
forms. We can assume that the path $\hat{h}_{k,t}$ is stationary
in $[0,\b]$ and $[1-\b,1]$. Again $\partial \hat{h}_{k,t}$ is
$\e/c_3$-transverse to zero and $\hat{h}_{k,t}$ has $C^3$-bounds
$c_3 c$ on $B_{\C^n}(0,1)$, where $c_3>0$ only depends on the
geometry of the manifold. On the other hand, note that
$\hat{h}_{k,t}$ is naturally defined all over $\C^n$. Now we
define the following section to the trivialized bundle
$\Phi_{t}^* \left( L^{\otimes k} \oplus L^{\otimes k} \right)$
on $B_{\C^n}(0,1)$:
 $$
 \hat{l}_{k,t}= (1,\l_t+ \hat{h}_{k,t}).
 $$

The goal now is to go back to the manifold. First, note that by
Lemma \ref{lem:asymp_Darboux} the chart $\Phi_{t}$ may be extended
as
 $$
 \Phi_{t}: B_{\C^n}(0, k^{1/2})  \to V_t,
 $$
where $V_t$ is a neighborhood of $p_t$ of $g_k$-radius
$O(k^{1/2})$, $p_t\in U_t \subset V_t$. Also
$s_{k,p_t}^{\text{\scriptsize ref}}$ is supported in a ball of
$g_k$-radius $O(k^{1/6})$.  Then
 $$
 u_{k,t}=  \left\{ \begin{array}{ll}
   {\displaystyle  \frac{\b-t}{\b} s_{k} +
 \frac{t}{\b} (\hat{l}_{k,0}\circ \Phi_{0}^{-1})\ox
  s_{k,p_0}^{\text{\scriptsize ref}}} , \qquad & t\in [0, \b] \\
 {\displaystyle (\hat{l}_{k,t}\circ \Phi_{t}^{-1})\ox
  s_{k,p_t}^{\text{\scriptsize ref}}} , \qquad &
  t\in [\b, 1-\b] \\
  {\displaystyle  \frac{t-1+\b}{\b} s_{k} +
 \frac{1-t}{\b} (\hat{l}_{k,1}\circ \Phi_{1}^{-1})\ox
  s_{k,p_1}^{\text{\scriptsize ref}}} , \qquad & t\in [1-\b, 1]
 \end{array} \right.
 $$
is a well-defined section of $L^{\ox k} \oplus L^{\ox k}$.
It is easy to check that
$u_{k,t}$ has $C^3$-bounds $c_4 c$, where $c_4$ depends only on
the geometry of the manifold, using that $|\bbd \Phi_t(z)| =
O(k^{-1/2}|z|)$ and that $s_{k,p_t}^{\text{\scriptsize ref}}$ has
Gaussian decay. Now in a small neighborhood of $p_t$ the map
 $$
 \p_{k,t}=\frac{u_{k,t}^2}{u_{k,t}^1}=\left\{ \begin{array}{ll}
 {\displaystyle \frac{(1-\mu) s^2_k + \mu\,
 s_{k,p_0}^{\text{\scriptsize ref}}
 (\l_0+ \hat{h}_{k,0}\circ \Phi_{0}^{-1})}{(1-\mu) s^1_k +
 \mu\, s_{k,p_0}^{\text{\scriptsize ref}}}}, \qquad &
 \mu=\frac{t}{\b}\in [0,1] \\
 \l_t + \hat{h}_{k,t} \circ \Phi_{t}^{-1}, &
 t \in [\b,1-\b] \\
 {\displaystyle \frac{(1-\mu) s^2_k + \mu\, s_{k,p_1}^{\text{\scriptsize ref}}
 (\l_1+ \hat{h}_{k,1}\circ \Phi_{1}^{-1})}{(1-\mu) s^1_k +
 \mu\, s_{k,p_1}^{\text{\scriptsize ref}}}}, & \mu=\frac{1-t}{\b}\in [0,1] \\
 \end{array} \right.
 $$
is well-defined, $\bd \p_{k,t}$ has a zero in $p_t$ and is
$\e/c_4$-transverse to zero in that small neighborhood. This is
clear in the second interval. In the first and third intervals it
is easy to check if we have previously normalized
$s_{k,p_t}^{\text{\scriptsize ref}}$ to be a positive real
multiple of $s_k^1$ at the point $p_t$. Note that $s_k^1$ is
uniformly bounded below since we are working off a neighborhood of
$Z(s_k^1)$.

Now we apply Proposition \ref{prop:get_tr} to the family $u_{k,t}$
choosing $\d=\min\{ \frac{\e}{4c_4}, \frac{\e}{2}\}$. We obtain a
new family $\hat{u}_{k,t}$ which satisfies the required
transversality conditions with a constant $\e'$ depending on $\d$.
The map $\hat{\p}_{k,t}=\frac{\hat{u}_{k,t}^2}{\hat{u}_{k,t}^1}$
is well defined outside $Z(\hat{u}_{k,t})$. It is clear that there
is a path of zeroes of $\bd \hat{\p}_{k,t}$ which is $C^0$-close
to $p_t$. Moreover, choosing $\a>0$ small, we define the family
 $$
 \tilde{s}_{k,t}= \left\{ \begin{array}{ll}
  {\displaystyle \frac{\a-t}{\a}s_{k,0} + \frac{t}{\a}\hat{u}_{k,0}},
   & t\in [0, \a] \\
 \hat{u}_{k,\frac{t-\a}{1-2\a}}, & t\in [\a,1-\a] \\
 {\displaystyle \frac{t-1+\a}{\a}s_{k,1}+\frac{1-t}{\a}\hat{u}_{k,1}},
 \qquad &
 t\in [1-\a, 1] \end{array} \right.
 $$
This defines a family of sections that are $\min(\e',
\e)$-transverse to zero in the sense that they satisfy the various
properties of Proposition \ref{prop:get_tr}.

Therefore the points $p_0$ and $p_1$ are joined by a path of
zeroes of $\partial \tilde{\p}_{k,t}$, where $\tilde{\p}_{k,t}$ is
the map associated to $\tilde{s}_{k,t}$. After applying
Proposition \ref{prop:get_model}, this path becomes a path of
critical points of the family of Lefschetz pencils. By the same
arguments as in the algebraic case the vanishing spheres
associated to the two critical points are conjugated through the
action of an orientation-preserving diffeomorphism.
The last point is to recall that the $C^3$-bounds and the
estimated transversality constants obtained along the way are
independent of the chosen points $p_0$ and $p_1$. This implies
that the constructions begin to work for a given $k$ for all
possible pairs of points. This concludes the proof.
\end{proof}

\begin{remark}
The above argument may be continued to prove that the Lagrangian
vanishing spheres are conjugated by symplectomorphisms of the
fiber. To check this, it suffices to use the parallel transport
defined by the symplectic orthogonal horizontal space in the
family of Lefschetz pencils on $M$. We skip this refinement to
prove directly the stronger statement of Theorem \ref{moser}.
\end{remark}

\section{Topology of the space of $\e$-transverse sections}
\label{comp_dual}

In the case of a complex projective manifold, where Lefschetz pencils
arise from lines in the dual projective space, the embedding
into it of the dual variety
has a further property that facilitates the study of its complement:

\begin{theorem} [Lefschetz hyperplane section theorem, Hamm--L\^e,
see \cite{Morse}] \label{t:hypsection}
 Let $M
\hookrightarrow \CP^N$ be a complex projective manifold, $M^* \hookrightarrow
(\CP^N)^*$ its dual variety,
and let $H \subset (\CP^N)^*$ be a generic linear subvariety of
complex dimension $d$. Then the homotopy group morphisms
$\pi_i(H-(H \cap M^*)) \to \pi_i((\CP^N)^*-M^*)$ induced by the
natural inclusion of spaces are isomorphisms for $i<d$, and an
epimorphism for $i=d$.
\end{theorem}

We are specially interested in the symplectic version of a
particular case of Theorem \ref{t:hypsection} in this paper:
whether the fundamental group of the space of regular values of a
Lefschetz pencil generates the fundamental group of the complement
of the dual variety. Our Theorem \ref{thm:topo} is somewhat weaker
than that, due to difficulties arising from the nature of the
symplectic dual variety (see Question \ref{question}).

\bigskip

Let $M$ be a symplectic manifold of integer class and
let $L\to M$ be the associated line bundle. For a compatible
almost complex structure $J$, $\e>0$, and $k$ large enough, define
\begin{equation}
  M^*_{\e,k,J}= \{ s:M \to L^{\ox k}: \text{$s$ is $\e$-transverse,
   $s$ has $C^3$-bounds $1$ (for $J$)} \}, \label{def:sect}
\end{equation}
with the natural $C^3$-topology. Observe that
$M^*_{\e_1,k,j}\subset M^*_{\e_2,k,J}$, for $\e_1\geq \e_2$. Note
first that the isotopy result of \cite{Au97} shows that for $k$
large enough (how large depending on $\e$), the zero sets of the
sections $s_k\in M^*_{\e,k,J}$ have isotopic zero sets
$W_k=Z(s_k)$.

Now we aim to understand better the fundamental group of the
spaces $M_{\e,k, J}^*$. We will study a suitably modified version
of it. Let $\G_s(L^{\ox k})$ denote the subspace of sections of
$L^{\ox k}$ whose zero locus is a smooth symplectic submanifold of
$M$, and let $\PP\G_s(L^{\ox k})$ be its projectivization.
Therefore the space $M_{\e,k,J}^*$ is an open subset of the space
$\G_s(L^{\ox k})$. We denote by $\rho_{\e,k,J}: M_{\e,k,J}^* \to
\PP\G_s(L^{\ox k})$ the natural map and set
 \begin{equation}\label{eqn:inva}
 \pi_{\e,k,J}=(\rho_{\e,k,J})_* (\pi_1(M_{\e,k,J}^*)) \subset \pi_1
 (\PP\G_s(L^{\ox k})).
 \end{equation}
Denote $\Map(W_k)=\text{Diff}^+(W_k)/ \text{Diff}_0(W_k)$ the
mapping class group of the fiber, i.e. its oriented diffeomorphisms modulo
isotopies. As in the algebraic case there is a geometric monodromy map
\begin{equation} \label{eq:geomono}
\mu : \pi_1(M_{\e,k,J}^*,*) \longrightarrow \Map(W_k)
\end{equation}
which factors through $(\rho_{\e,k,J})_*$, thus defining a map
$\pi_{\e,k,J} \to \Map(W_k)$.

The proposed way for studying the topology of $M_{\e,k,
J}^*$ is by means of a Lefschetz hyperplane theorem that allows a
finite dimensional reduction. The idea is to restrict ourselves to
a linear subspace of finite dimension. More specifically, choose
sufficiently generic sections $s_k^0,\ldots, s_k^N$ in
$M^*_{\e,k,J}$ and let $V_{k,N}$ be their linear span. Then one
may ask whether there is an isomorphism $\pi_j(V_{k,N} \cap
M^*_{\e,k,J})\iso \pi_j(M^*_{\e,k,J})$ for $j$ below the middle
dimension as in Theorem \ref{t:hypsection}. In
the case of the fundamental group we consider a symplectic
Lefschetz pencil $\p_k=s_k^2/s_k^1$, which defines a natural map:
\begin{eqnarray*}
  \Psi_k:\CP^1- \{p_1, \ldots, p_l \} & \to & \PP\G_s(L^{\ox k}) \\
   \left[ q_1 , q_2 \right] & \mapsto & [q_1s_k^2-q_2s_k^1]
\end{eqnarray*}
where $p_1, \ldots, p_l\in \CP^1$ are the images of the critical
points.

\begin{theorem} \label{thm:topo}
 Let $\e>0$ be small enough.
 Then there exists $k_0$ such that for any $k \geq k_0$
 the image of the geometric monodromy map
 $$
  \mu: \pi_1(M^*_{\e,k,J})\to \Map(W_k)
 $$
 lies in the subgroup of $\Map(W_k)$ generated by the
 positive generalized Dehn twists of a single
 Donaldson's Lefschetz pencil $\phi_k :
 M - N_k \to \CP^1$.
\end{theorem}

\begin{proof}
First we need to introduce some definitions. Abusing notation we
will denote by $L$ the pull-back of the bundle $L$ to the manifold
$M\x S^1$. Then we can define $C^3$-bounds of a sequence of
sections of $L^{\ox k}\to M\x S^1$ just by requiring that the
restriction to each fiber of $\pi:M\x S^1 \to S^1$ satisfy these
bounds. In the same way the section $s:M\x S^1 \to L^{\ox k}$ is
$\e$-transverse to zero if the restriction to the fiber
$\pi^{-1}(z)$ is $\e$-transverse for any $z\in S^1$.

We fix, once for all, as a base point, a section $s_{k,0}^2 \in
M^*_{\e,k,J}$, that exists by \cite{Do96} for some small $\e>0$.
Take an element $A\in \pi_1(M^*_{\e,k,J}, s_{k,0}^2)$. Fix a
representative $\a:S^1 \to M^*_{\e,k,J}$, i.e.,\ $[\a]=A$. The map
$\a$ defines canonically a section $s_k^2$ of $L^{\ox k} \to M\x
S^1$ by the formula
 $$
 s_k^2(x,e^{2 \pi i t})=(\a(e^{2\pi i t}))(x).
 $$
with the base point condition that $s_k^2(\cdot\, ,1)=s_{k,0}^2$.
It is clear that $s_k^2$ has $C^3$-bounds $1$ and it is
$\e$-transverse to zero. Now, choose another $\e$-transverse
section $s_k^1$ of $L^{\ox k} \to M$ such that $s_{k}^1 \oplus
s_{k,t}^2$ is a section of  $L^{\ox k} \oplus L^{\ox k}$ with
$C^3$-bounds $1$ for all $t \in S^1$. Moreover we may also
assume that $s_k^1\oplus s_{k,j}^2$ is already an $\e$-transverse
Lefschetz pencil (for $\e$ small enough), for $j=0,1$. Denote
again by $s_k^1$ the pull-back of the section to $M\x S^1$.
Applying the $S^1$-parametric version of Propositions
\ref{prop:get_tr} and \ref{prop:get_model} we perturb $s_k^1\oplus
s_k^2$ into $\s_k^1 \oplus \s_k^2$, a new section that defines a
$S^1$-family of $\e'$-transverse Lefschetz pencils, for some
$\e'<\e$.

Note that in the process of Proposition \ref{prop:get_tr} it is
necessary to perturb only $s_{k,t}^2$. This is because the
condition (ii) is already satisfied by $s_{k,t}^1=s_k^1$, the
condition (iv) is obtained in \cite{Do99} only perturbing
$s_{k,t}^2$ and the proof of the condition (iii) can be changed
using the alternative proof written in \cite{Au97} that does not
need to perturb $s_{k}^1$. So $\s_{k}^1=s_k^1$. Moreover, the
perturbation in $s_{k,t}^2$ can be very small so that the sequence
of sections $\s_{k,t}^2$ is isotopic to the previous one. Also we
may assume that $\s_{k,j}^2=s_{k,j}^2$ for $j=0,1$, since
$s_k^1\oplus \s_{k,j}^2$ are already Lefschetz pencils.

Now we want to compute the monodromy associated to the zero set of
$\s_k^2$, which is the same as the one of $s_k^2$.

\begin{figure}[h]
\begin{center}
\includegraphics{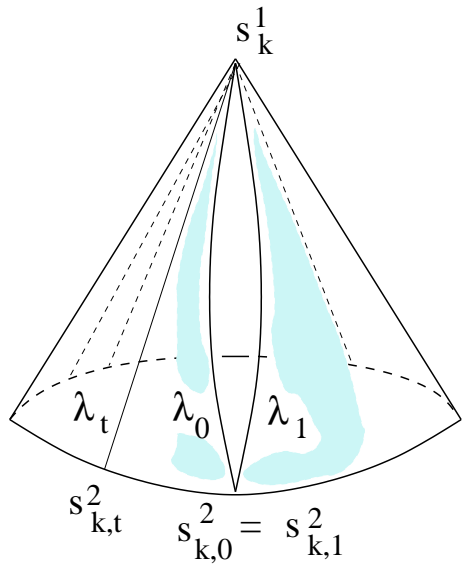}
\end{center}
\caption{Homotopy between $\s^2_{k,t}$ and $\lambda_0^{-1}
\cdot \lambda_1$ in $M^*_{\e,k,J}$}
\end{figure}

For each fiber $M\x \{e^{2 \pi i t}\}$ we know that the images of
the critical points of the associated Lefschetz pencil
$\phi_{k,t}$ are a finite number of disjoint points. For the whole
family of pencils over $S^1$ the critical values describe a braid
in $\CP^1$, so we can choose a continuous family of paths $\l_t :
[0,1] \to \CP^1$ joining $0$ with $\infty$, such that each $\l_t$
avoids the critical values of the pencil $\phi_{k,t}$.

Now, the disk bounded by $\s_{k,t}^2$, $\l_0$, $\l_1$ and the
constant path $s_k^1$ defines a base-pointed homotopy between the
loops $\s_{k,t}^2$ and $\l_0^{-1} \cdot \l_1$. The latter loop is
contained in the regular values of the Lefschetz pencil
$\phi_{k,0}= [s_k^1,\s_{k,0}^2]$, therefore the monodromy of the
family of sections of $M$ over it is a product of direct
generalized Dehn twists.

As $\e'$ depends on $\e$, we may find $k_\e$ only dependent on
$\e$ such that, for any $k\geq k_\e$, all $\e'$-transverse pencils
are isotopic. Furthermore the isotopy between two pencils $[s_k^1,
s_{k,0}^2]$ and $[\s_k^1,s_{k,0}^2]$ of $M$ may be achieved by a
family of pencils that keeps fixed the section $s_{k,0}^2$. Thus
the image of the geometric monodromy map $\mu$ is contained in the
subgroup of $\Map (W_k)$ generated by a single pencil, which may
be assumed to be $\e$-transverse.
\end{proof}

\section{Nontriviality of the vanishing spheres} \label{case:moser}

Now we will prove Theorems \ref{moser} and \ref{non-triv}. The
first one is an immediate consequence of Theorem \ref{thm:topo}.

\noindent {\it Proof of Theorem \ref{moser}.} Recalling the proof
of Theorem \ref{main:alg}, take two vanishing spheres $S_0$, $S_1$
in a fiber of a Donaldson Lefschetz pencil and join them by a path
of Lagrangian spheres $S_t$ such that $S_t \subset W_t'$, where
$W_t'$ is a family of fibers in a family of pencils $\phi_{k,t}$.
This means that $W_t'=Z(s_{k,t})$, $s_{k,t} \in M_{\e,k,J}^*$. We
can suppose that $s_{k,0}=s_{k,1}$ is the base point. Now
$s_{k,t}$ can be homotoped, by Theorem \ref{thm:topo}, to a path
inside a Lefschetz pencil whose zero set we denote by $W_t$. As
usual, parallel transport along the path generates a family of
symplectomorphisms $\phi_{t,s}: W_t \to W_{t+s}$. The
symplectomorphism $\phi_{0,1}$ is generated by the composition of
the generalized Dehn twists associated to the critical points of
the pencil bounded by the path.

However $S_0$ is not, in general, preserved by $\phi_{0,1}$.
Denote by $S_{t,s}=\phi_{t,s}(S_t)$. Then we have that the
continuous family $S_{t,1-t}$ provides a Lagrangian isotopy in
$W_0=W_1$ between $S_1=S_{1,0}$ and $S_{0,1}$. This shows that the
the spheres $S_0$ and $S_1$ are conjugated by $\phi_{0,1}$, up to
Lagrangian isotopy. \hfill $\Box$

\begin{remark} \label{rem:isotopy}
The spheres are conjugated only up to Lagrangian isotopy. We can
also claim that they are conjugated up to Hamiltonian isotopy
because the obstruction of a Lagrangian isotopy to be lifted to a
Hamiltonian isotopy lies in the family of groups $H^1(L_t)$ (cf.\
Exercise 6.1.A in \cite{Po01}). This is an infinitesimal condition
to be fulfilled for each member of the family. In our case this
cohomology groups vanish in the case of dimension strictly bigger
than $4$. In dimension $4$, we can check that the spheres at the
starting and end points of the family are exact (using the local
models close to critical points) and this is enough to get a
Hamiltonian isotopy (recall that in dimension $4$ the fibers are
Riemann surfaces).
\end{remark}

Theorem \ref{moser} admits an extension from the symplectomorphism
group of the fiber to that of the total space.

\begin{theorem} \label{thm:sympglob}
For a Donaldson's pencil on $M$ with $k$ large enough, fix a
generic fiber $W$. If $S_0,S_1 \subset W$ are Lagrangian spheres
that appear as vanishing cycles of the pencil then there exists a
symplectomorphism $\Phi : M \rightarrow M$ such that $\Phi(W)=W$
and $\Phi(S_0)$ equals $S_1$ up to Lagrangian isotopy in the fiber
$W$.
\end{theorem}

\begin{proof}
It follows from Proposition \ref{prop:hard} that there exists a
loop of approximately holomorphic sections $s_{k,t}$ such that for
$t=0,1$ the zero set $Z(s_{k,t})=W$, and there is a family of
Lagrangian spheres $S_t$ in $W_t=Z(s_{k,t})$ connecting the two
given vanishing spheres.

Now by Proposition 4 in \cite{Au97}, there is a continuous family
of symplectomorphisms $\phi_t:M\to M$ such that $\phi_t(W_0)=W_t$.
Again, we cannot assume that $S_t$ is preserved by the family.
However, if we construct $S_{t,s}=\phi_s(S_t)\subset W_{t+s}$, we
have that the family $S_{t,1-t}$ defines a Lagrangian isotopy in
$W_1=W_0$ between $S_1$ and $S_{0,1}=\phi_1(S_0)$. So, we have
shown, up to Lagrangian isotopy in the fiber, that $\phi_1$ is the
required symplectic map.
\end{proof}

\noindent {\em Proof of Theorem \ref{non-triv}.\/} The proof of
Theorem \ref{non-triv} works looking for a contradiction. Theorem
\ref{moser} implies that all the vanishing spheres are either
homologically trivial or homologically non-trivial. Let us suppose
that they are homologically trivial.

Let $\p_k:M-N_k\to \CP^1$ be the Lefschetz pencil obtained from
Theorem \ref{thm:exist} for $k$ large, and denote by $F_k$ a regular
fiber. We blow-up the manifold $M$
along $N_k$ to obtain a fibration $\tilde{\p}_k: \tilde{M}_k \to
\CP^1$.

Now let us compute the (growth of the) Betti numbers of
$\tilde{M}_k$. First, by \cite[Proposition 5]{Au97}, $\chi(F_k)=
(-1)^{n-1} \vol(M) k^n+ O(k^{n-1})$, where $\vol(M)=\int_M \left(
\frac{\omega}{2\pi} \right)^n$. Now the Lefschetz hyperplane
theorem \cite{Au97} implies that $b_i(F_k) = b_{i}(M)$, which is
uniformly bounded for $i<n-1$ (i.e.,\ independently of $k$). By
Poincar{\'e} duality, $b_{2n-2-i}(F_k)$ is also uniformly bounded.
Therefore
  \begin{equation}\label{eqn:p2}
  b_{n-1}(F_k)=\vol(M) k^n+ O(k^{n-1}).
  \end{equation}
A handlebody decomposition of $\tilde{M}_k$ is obtained as
follows: take a tubular neighbourhood $\nu(F_k)$ of a smooth fiber
$F_k$ of $\tilde{\p}_k$, and attach $n$-handles along
$(n-1)$-spheres embedded in the boundary  $\bd \nu(F_k)=F_k \times
S^1$ as the vanishing cycles at different fibers. Since all of
these spheres are homologically trivial, the homology group
$H_{n-1}$ does not change. So $H_{n-1}(\tilde{M}_k-F_k\times
D^2)\iso H_{n-1}(F_k)$. Finally, we have to attach $F_k \times
D^2$ to get $\tilde{M}_k$. We have a Mayer-Vietoris exact sequence
  $$
   H_{n-1}(F_k\x S^1) \to H_{n-1}(\tilde{M}_k-F_k\times
D^2) \oplus H_{n-1}(F_k\x D^2) \to H_{n-1}(\tilde{M}_k) \to
H_{n-2} (F_k\times S^1)
  $$
Using that $b_{n-2}(F_k \times S^1)$ is uniformly bounded, that
$H_{n-1}(F_k\x S^1) \to H_{n-1}(F_k\x D^2)$ is surjective with
kernel uniformly bounded and \eqref{eqn:p2}, we get
  \begin{equation}\label{eqn:p}
  b_{n-1}(\tilde{M}_k ) = \vol(M) k^n+ O(k^{n-1}).
  \end{equation}

Analogously to the computation of $F_k$, and using that
$\chi(N_k)= (-1)^{n-2} (n-1)\vol(M) k^n +O(k^{n-1})$, we get that
all the Betti numbers of $N_k$ are bounded except for
$b_{n-2}(N_k)=(n-1)\vol(M) k^n+ O(k^{n-1})$. In the blow up
$\tilde M_k$, the exceptional divisor $\tilde{N}_k$ is a fibration
over $N_k$ by $\CP^1$. There is a spectral sequence with
$E_2$-term $H_*(N_k)\ox H_*(S^2)$ and abutting to
$H_*(\tilde{N}_k)$. This gives that all $b_i(\tilde{N}_k)$ are
bounded, for $i\neq n-2, n$, and that
$b_{n-2}(\tilde{N}_k)=b_n(\tilde{N}_k)=(n-1)\vol(M) k^n+
O(k^{n-1})$.

Consider the following pieces of long exact sequences in homology,
$$
\begin{array}{ccccccc}
  H_{n-1}(\tilde{N}_k) & \stackrel{\tilde{i}_*}{\too}  &
  H_{n-1}(\tilde{M}_k) & \stackrel{\tilde{j}_*}{\too} &
  H_{n-1}(\tilde{M}_k, \tilde{N}_k) & \stackrel{\tilde{\bd}_*}{\too}
  &  H_{n-2}(\tilde{N}_k)  \\
  \downarrow && \downarrow && \downarrow && \downarrow \\
  H_{n-1}({N}_k) & \stackrel{{i}_*}{\too}   &   H_{n-1}(M) & \stackrel{j_*}{\too} &
  H_{n-1}(M, {N}_k) & \stackrel{\bd_*}{\too} & H_{n-2}({N}_k)
\end{array}
$$

The kernel of $\tilde{i}_*$ is uniformly bounded, since
$b_{n-1}(\tilde{N}_k)$ is so. Also $\im \tilde{j}_* =\ker
\tilde{\bd}_* \subset \ker \bd_* = \im j_*$ is uniformly bounded.
Therefore, $b_{n-1}(\tilde{M}_k)$ is uniformly bounded. This is a
contradiction with \eqref{eqn:p}. \hfill $\Box$

\vskip5pt

\section{Towards a symplectic invariant} \label{invariant}

In the case of complex projective manifolds, the topology of the
space of smooth sections of the ample bundle $L \rightarrow M$ is
governed by the hyperplane section theorem of Hamm--L\^e (Theorem
\ref{t:hypsection}), which in particular shows how the monodromy
of generic families of dimension $1$ or $2$ of hyperplane sections
capture the geometric monodromy of all families of sections.

In order to extend these properties from the complex algebraic to the
symplectic case and define a new set of symplectic invariants, the
next natural step is to prove that the groups $\pi_{\e,k,J}$ defined in
\eqref{eqn:inva}, are independent
of $J$ and of $\e$ for $k$ large enough. The authors have not been
able to prove this, thus can only pose the following

\begin{question} [symplectic Hamm--L\^e hyperplane sections]
\label{question} Does the morphism $\pi_l(\CP^j \cap
\PP\G_s(L^{\ox k})) \rightarrow \pi_l(\PP \Gamma_s(L^{\ox k}))$
induce an isomorphism $\pi_l(\CP^j \cap \PP\G_s(L^{\ox k}))
\rightarrow (\rho_{\e,k,J})_*\pi_l(M_{\e,k,J}^*)$ for $l<j$ and an
epimorphism for $l=j$, for a $(j+1)$-dimensional ``generic''
family of Donaldson's $\e$-transverse sections?
\end{question}

An affirmative answer would yield a family of invariants, the
$\e$--transverse homotopy groups of the complement $M_{\e,k,J}^*$
of the dual variety, which would be an extension of the Auroux and
Katzarkov invariants of \cite{AK00}. It must be warned that the
only evidence available for an affirmative answer is the case of
complex algebraic manifolds, and the authors have been able to
make only a modest extension to the general symplectic case:

\begin{corollary}
For $\e>0$ small and $k$ large enough, we have $(\rho_{\e,k,J})_*
(H_1(M_{\e,k,J}^*)) =\Z/ m\Z \subset H_1 (\PP\G_s(L^{\ox k}))$,
where $m$ is a divisor of the number $n_k$ of singular fibers in
the Lefschetz pencil $\phi_k : M-N_k \to \CP^1$.
\end{corollary}

\begin{proof}
By Theorem \ref{thm:topo} the image of $H_1(M_{\e,k,J}^*)$ is
generated by some compositions of monodromies around critical
points of a single Donaldson's Lefschetz pencil. But the proof of
Proposition \ref{prop:hard} shows that all the loops around
different critical points are freely homotopic. Therefore this
abelian group is cyclic, generated by the $1$-cycle $\gamma$
bounding a small disk $D$ around a critical value of the pencil.

The inclusion of the regular values of any Donaldson's Lefschetz
pencil in $\PP\G_s(L^{\ox k})$ induces a homology $n_k \gamma \sim
0$.
\end{proof}

The ``natural'' value for $m$ above should be the number $n_k$.

\end{document}